\documentclass[11pt]{article}
\usepackage{amsmath}
\usepackage{amsfonts}
\usepackage{amsthm}
\usepackage{verbatim}
\begin{document}
\def\N{\mathbb{N}}
\def\F{\mathbb{F}}
\def\Z{\mathbb{Z}}
\def\R{\mathbb{R}}
\def\Q{\mathbb{Q}}
\def\H{\mathcal{H}}
\parindent= 3.em 
\parskip=5pt

\centerline{\bf{ ON A QUARTIC EQUATION AND TWO FAMILIES}}
\centerline{\bf{ OF HYPERQUADRATIC CONTINUED FRACTIONS }}
\centerline{\bf{ IN POWER SERIES FIELDS}}
\centerline{\bf{by}}\centerline{\bf{Kh. Ayadi and A. Lasjaunias}}
\noindent{\bf{1. Introduction}} \par Throughout this note $p$ is an odd prime number and $\F_p$ is the finite field with $p$ elements. We consider an indeterminate $T$ and $\F_p((T^{-1}))$, the field of power series in $1/T$ over the finite field $\F_p$, here simply denoted by $\F(p)$. A non-zero element of $\F(p)$ is $$\alpha=\sum_{i\leq i_0}u_iT^i \quad \text{ where }\quad i\in \Z, u_i\in \F_p \quad \text{ and }\quad u_{i_0}\neq 0.$$An ultrametric absolute value is defined over this field by $\mid \alpha \mid =\mid T \mid^{i_0}$ where $\mid T\mid$ is a fixed real number greater than 1. We will also consider the subset $\F(p)^+=\lbrace \alpha \in \F(p) \quad s.t. \mid \alpha \mid >1\rbrace$. Note that $\F(p)$ is the completion of the field $\F_p(T)$ for this absolute valute. The fields $\F(p)$ are analogues of the field of real numbers, consequently many questions in number theory in the context of real numbers, such as Diophantine approximation and continued fractions, can be transposed in the frame of formal power series which is considered here. We are concerned with continued fractions for elements of this field $\F(p)$ which are algebraic over the field $\F_p(T)$. 
\par The starting point of our work is a particular quartic equation, with coefficients in $\F_p[T]$, where $p$ is an arbitrary prime greater than 3. This algebraic equation is the following:
$$(9/32)X^4-TX^3+X^2-8/27=0.\eqno{(Eq)}$$The origin of this equation is due to Mills and Robbins [MR]. Mills and Robbins actually considered another equation: $(Eq_1)\quad X^4+X^2-TX+1=0$. The very simple form of this last equation explains why it was considered by chance, while searching for promising algebraic continued fraction expansions. Using a computer, Mills and Robbins observed that $(Eq_1)$ has a root in $\F(13)$ presenting a remarkable continued fraction expansion. This continued fraction expansion could only be partially conjectured in [MR] and only later, in a complicated form, fully conjectured in [BR]. Finally, the conjecture concerning the continued fraction of the solution of $(Eq_1)$ in $\F(13)$ was proved in [L3]. In [L5], it has been remarked that Mills and Robbins equation should be considered in all chracteristic $p>3$, by reading $X^4+X^2-TX-1/12=0$. After a transformation, this led to the above equation $(Eq)$ (see [L5, p. 30]). 
\par For each prime $p>3$, $(Eq)$ has a unique root in $\F(p)^+$, denoted by $\alpha(p)$. This root can be expanded as an infinite continued fraction. The continued fraction for $\alpha(p)$ varies according to the value of $p$ but, for all $p$, it appears to have a singular pattern. Moreover, observations by computer show that there are two different general patterns, according to the case considered: $p\equiv 1 \mod 3$ or $p\equiv 2 \mod 3$. The first case, $p\equiv 1 \mod 3$ (and particularly $p=13$), has been extensively studied by the second author and this study has generated different works in the area of continued fractions in power series fields  ([L2], [L3], [L4] and [L5]). To show the differences and the similarities between both cases, we will recall several results already known for the first case. It will appear that the continued fraction expansion for the solution of $(Eq)$ belongs to two large families of expansion, according to the remainder of $p$ modulo $3$. Hence, the great interest of our equation will be to give us the opportunity to introduce and to describe these families. 
\par In order to illustrate our subject and to provoke the curiosity of the reader, we show, at the end of this introduction, what could be seen on a computer screen when considering the first few hundreds of partial quotients of the solution $\alpha(p)$ of $(Eq)$, for the first two values of $p$. Note that the partial quotients, which are quickly very large, are only represented there by their degree in the indeterminate $T$. See Figure~\ref{solp5} and Figure~\ref{solp7}. 

\begin{figure}[b]
\caption{C.f.e. for the solution $\alpha(5)$ of $(Eq)$ (450 p. q.)\label{solp5}}
\small
\begin{verbatim}
[1, 1, 1, 1, 1, 1, 1, 1, 1, 1, 1, 1, 9, 41, 9, 1, 1, 1, 1, 1, 1, 1, 1, 
1, 1, 1, 1, 1, 1, 1, 1, 9, 1, 1, 1, 1, 1, 1, 1, 1, 1, 1, 1, 1, 1, 1, 1, 
1, 9, 1, 1, 1, 1, 1, 1, 1, 1, 1, 1, 1, 1, 1, 1, 1, 1, 9, 41, 9, 1, 1, 1, 
1, 1, 1, 1, 1, 1, 1, 1, 1, 1, 1, 1, 1, 9, 1, 1, 1, 1, 1, 1, 1, 1, 1, 1,
1, 1, 1, 1, 1, 1, 9, 1, 1, 1, 1, 1, 1, 1, 1, 1, 1, 1, 1, 1, 1, 1, 1, 9,
41, 209, 1041, 209, 41, 9, 1, 1, 1, 1, 1, 1, 1, 1, 1, 1, 1, 1, 1, 1, 1,
1, 9, 1, 1, 1, 1, 1, 1, 1, 1, 1, 1, 1, 1, 1, 1, 1, 1, 9, 1, 1, 1, 1, 1,
1, 1, 1, 1, 1, 1, 1, 1, 1, 1, 1, 9, 41, 9, 1, 1, 1, 1, 1, 1, 1, 1, 1, 1,
1, 1, 1, 1, 1, 1, 9, 1, 1, 1, 1, 1, 1, 1, 1, 1, 1, 1, 1, 1, 1, 1, 1, 9,
1, 1, 1, 1, 1, 1, 1, 1, 1, 1, 1, 1, 1, 1, 1, 1, 9, 41, 9, 1, 1, 1, 1, 1,
1, 1, 1, 1, 1, 1, 1, 1, 1, 1, 1, 9, 1, 1, 1, 1, 1, 1, 1, 1, 1, 1, 1, 1,
1, 1, 1, 1, 9, 1, 1, 1, 1, 1, 1, 1, 1, 1, 1, 1, 1, 1, 1, 1, 1, 9, 41,
209, 41, 9, 1, 1, 1, 1, 1, 1, 1, 1, 1, 1, 1, 1, 1, 1, 1, 1, 9, 1, 1, 1,
1, 1, 1, 1, 1, 1, 1, 1, 1, 1, 1, 1, 1, 9, 1, 1, 1, 1, 1, 1, 1, 1, 1, 1,
1, 1, 1, 1, 1, 1, 9, 41, 9, 1, 1, 1, 1, 1, 1, 1, 1, 1, 1, 1, 1, 1, 1, 1,
1, 9, 1, 1, 1, 1, 1, 1, 1, 1, 1, 1, 1, 1, 1, 1, 1, 1, 9, 1, 1, 1, 1, 1,
1, 1, 1, 1, 1, 1, 1, 1, 1, 1, 1, 9, 41, 9, 1, 1, 1, 1, 1, 1, 1, 1, 1, 1,
1, 1, 1, 1, 1, 1, 9, 1, 1, 1, 1, 1, 1, 1, 1, 1, 1, 1, 1, 1, 1, 1, 1, 9,
1, 1, 1, 1, 1, 1, 1, 1, 1, 1, 1, 1, 1, 1, 1, 1, 9, 41, 209, 41, 9, 1, 1]
\end{verbatim}
\end{figure}

\begin{figure}[b]
\caption{C.f.e. for the solution $\alpha(7)$ of $(Eq)$ (430 p. q.)\label{solp7}}
\small
\begin{verbatim}
[1, 1, 1, 3, 1, 1, 1, 1, 3, 1, 1, 1, 1, 3, 1, 1, 1, 1, 17, 1, 1, 1, 1,
3, 1, 1, 1, 1, 3, 1, 1, 1, 1, 3, 1, 1, 1, 1, 3, 1, 1, 1, 1, 17, 1, 1, 1,
1, 3, 1, 1, 1, 1, 3, 1, 1, 1, 1, 3, 1, 1, 1, 1, 3, 1, 1, 1, 1, 17, 1, 1,
1, 1, 3, 1, 1, 1, 1, 3, 1, 1, 1, 1, 3, 1, 1, 1, 1, 3, 1, 1, 1, 1, 115,
1, 1, 1, 1, 3, 1, 1, 1, 1, 3, 1, 1, 1, 1, 3, 1, 1, 1, 1, 3, 1, 1, 1, 1,
17, 1, 1, 1, 1, 3, 1, 1, 1, 1, 3, 1, 1, 1, 1, 3, 1, 1, 1, 1, 3, 1, 1, 1,
1, 17, 1, 1, 1, 1, 3, 1, 1, 1, 1, 3, 1, 1, 1, 1, 3, 1, 1, 1, 1, 3, 1, 1,
1, 1, 17, 1, 1, 1, 1, 3, 1, 1, 1, 1, 3, 1, 1, 1, 1, 3, 1, 1, 1, 1, 3, 1,
1, 1, 1, 17, 1, 1, 1, 1, 3, 1, 1, 1, 1, 3, 1, 1, 1, 1, 3, 1, 1, 1, 1, 3,
1, 1, 1, 1, 115, 1, 1, 1, 1, 3, 1, 1, 1, 1, 3, 1, 1, 1, 1, 3, 1, 1, 1,
1, 3, 1, 1, 1, 1, 17, 1, 1, 1, 1, 3, 1, 1, 1, 1, 3, 1, 1, 1, 1, 3, 1, 1,
1, 1, 3, 1, 1, 1, 1, 17, 1, 1, 1, 1, 3, 1, 1, 1, 1, 3, 1, 1, 1, 1, 3, 1,
1, 1, 1, 3, 1, 1, 1, 1, 17, 1, 1, 1, 1, 3, 1, 1, 1, 1, 3, 1, 1, 1, 1, 3,
1, 1, 1, 1, 3, 1, 1, 1, 1, 17, 1, 1, 1, 1, 3, 1, 1, 1, 1, 3, 1, 1, 1, 1,
3, 1, 1, 1, 1, 3, 1, 1, 1, 1, 115, 1, 1, 1, 1, 3, 1, 1, 1, 1, 3, 1, 1,
1, 1, 3, 1, 1, 1, 1, 3, 1, 1, 1, 1, 17, 1, 1, 1, 1, 3, 1, 1, 1, 1, 3, 1,
1, 1, 1, 3, 1, 1, 1, 1, 3, 1, 1, 1, 1, 17, 1, 1, 1, 1, 3, 1, 1, 1, 1, 3,
1, 1, 1, 1, 3, 1, 1, 1, 1, 3, 1, 1, 1, 1, 17, 1, 1, 1, 1, 3, 1, 1, 1, 1]
\end{verbatim}}
\end{figure}

\par If the solution of $(Eq)$ has a peculiar continued fraction expansion, for each $p>3$, this is due to the fact that this element is hyperquadratic. Let $t\geq 0$ be an integer and $r=p^t$, an irrational element of $\F(p)$ will be called hyperquadratic of order $t$ if it satisfies a non-trivial algebraic equation of the following form$$uX^{r+1}+vX^r+wX+z=0 \quad \text{ where }\quad (u,v,w,z)\in (\F_p[T])^4.$$Note that a hyperquadratic element of order $0$ is simply irrational quadratic. We shall see that the solution of $(Eq)$ is hyperquadratic of order $1$ if $p\equiv 1 \mod 3$ and hyperquadratic of order $2$ if $p\equiv 2 \mod 3$. The reader may consult the introduction of [BL] for more precisions and references on hyperquadratic elements. Hyperquadratic power series in $\F(p)$ have long been considered by mathematicians studying Diophantine approximation in function fields of positive characteristic, such as Mahler [M], Osgood [O], Voloch [V] and de Mathan [dM]. Simultaneously, other mathematicians, such as Baum and Sweet [BS] or Mills and Robbins [MR], have observed that the continued fraction expansion of certain hyperquadratic elements could be explicitly given. For a survey on the different contributions of these mathematicians in this area, the reader is refered to [L1]. For a good account on continued fractions and Diophantine approximation in power series fields, aswell as more references, see Schmidt's article [Sc].
\par We shall now describe briefly the continued fraction for the solution of $(Eq)$ in $\F(p)^+$. This root is expanded as an infinite continued fraction $\alpha(p)=[a_1,a_2,\dots,a_n,\dots]$, where the partial quotients $a_i$ are non-constant polynomials in $\F_p[T]$. From the equation, one can obtain the begining of the power series expansion and we have $\alpha(p)=(32/9)T-1/T+\dots$, which implies that $a_1=(32/9)T$. To describe the sequence $(a_n)_{n\geq 1}$, in both cases, we need to introduce a particular polynomial in $\F_p[T]$ and two sequences of polynomials related to it.
\par Throughout this note $p$ is a prime with $p\geq 3$ (except when we consider $(Eq)$, where $p>3$) and $k$ an integer with $1\leq k<p/2$. Then we define $P_k(T)=(T^2-1)^k\in \F_p[T]$. From $P_k$, we introduce in $\F_p[T]$ two sequences of polynomials  $(A_{n})_{n\geq 0}$ and $(B_{n})_{n\geq 0}$ as follows. The first one is defined by$$A_{0}=T \quad \text{ and recursively }\quad  A_{n+1}=[A_{n}^p/P_k] \quad \text{ for }\quad n\geq 0.$$ Here the brackets denote the integral (i.e. polynomial) part of the rational function. While the second one is defined by
$$B_{0}=A_{0}=T \quad \text{ and }\quad B_{1}=A_{1}=[T^p/P_k]$$ 
and recursively $$B_{n+1}=B_{n}^pP_k^{(-1)^{n+1}}\quad \text{ for }\quad n\geq 1.$$ 
We are particularly interested in the degrees of these polynomials. We set $u_n=\deg(A_{n})$ and $v_n=\deg(B_{n})$. From the recursive definition of these polynomials, we get $u_0=v_0=1$ and also 
$$u_{n+1}=pu_n-2k \quad \text{ and }\quad v_{n+1}=pv_n+2k(-1)^{n+1}\quad \text{ for }\quad n\geq 0.$$ 
Note that the sequence $(u_n)_{n\geq 0}$ is constant if $2k=p-1$, then we have $A_n=T$ for $n\geq 0$. Otherwise, 
both sequences  $(u_n)_{n\geq 0}$ and $(v_n)_{n\geq 0}$ are strictly increasing. 
\newline Note that, for $p=5$ and $k=2$, we obtain : ${\bf{v}}=1,1,9,41,209,1041,\dots$ Whereas, for $p=7$ and $k=2$, we obtain ${\bf{u}}=1,3,17,115,801,\dots$
\newline (See below and also Figure 1 and Figure 2.)
\par In the first case, $p\equiv 1 \mod 3$, we set $k=(p-1)/3$ and we consider the sequence $(A_{n})_{n\geq 1}$, introduced above. In [L5], it has been proved (with a bound on the prime number $p$) that there exists a sequence $(\lambda_n)_{n\geq 1}$ in $\F_p^*$ and a sequence $(i(n))_{n\geq 1}$ in $\N$, such that
$$a_n=\lambda_nA_{i(n)}\quad \text{ for }\quad n\geq 1.\eqno{(I)}$$ 
Both sequences $(\lambda_n)_{n\geq 1}$ and $(i(n))_{n\geq 1}$ have been given explicitly (see [L4]).
\par In the second case, $p\equiv 2 \mod 3$, we set $k=(p+1)/3$ and we consider the sequence $(B_{n})_{n\geq 1}$ introduced above. Our observation, based on computer calculations giving a finite number of partial quotients, implies the following conjecture: there exists a sequence $(\lambda_n)_{n\geq 1}$ in $\F_p^*$ and a sequence $(i(n))_{n\geq 1}$ in $\N$, such that$$a_n=\lambda_nB_{i(n)}\quad \text{ for }\quad n\geq 1.\eqno{(II)}$$ 
\par In the first case, the formulas giving the sequence $(\lambda_n)_{n\geq 1}$ are quite sophisticated, as can be seen for instance for $p=13$ in [BR] and [L3]. Moreover, our proof and, consequently, the method to obtain this sequence are complicated (see [L2] and [L4]). For these reasons, in the second case, we have not tried to describe the sequence $(\lambda_n)_{n\geq 1}$. Moreover, still in this second case, we have decided to describe the sequence $(i(n))_{n\geq 1}$ as a consequence of a conjecture about more general continued fractions. The tools used to obtain a proof, in the first case, might well be applied in the second case, but we are aware that a different approach would be desirable. This note is complementary 
to [L5], and hopefully it may shed new light on this mysterious quartic equation. 
\par We sketch here the organization of this work. To obtain the proof in the first case, it has been necessary to consider hyperquadratic continued fractions more general than the one of the root of $(Eq)$. It happens that the same argument is true for the second case. In the next section we shall introduce these families, which we will call $P_k$-expansions of first kind and of second kind. In section 3, we will define and describe partially some $P_k$-expansions, which we call perfect. In section 4, we show that the solution of $(Eq)$ is a perfect $P_k$-expansion of first kind if $p\equiv 1 \mod 3$ and of second kind if $p\equiv 2 \mod 3$. In the last section, we give a mesure of the growth of the degrees of the partial quotients (the irrationality measure of the continued fraction) for the perfect $P_k$-expansions in both cases. We apply it to the solution of $(Eq)$ and we get the irrationality measure for $\alpha(p)$, equal to $8/3$ in the first case and equal to $4$ in the second one. 

\clearpage

\noindent{\bf{2. $P_k$-expansions }}
\par Concerning continued fractions in this area, we use classical notation, as they can be found for instance in the second section of [LY]. Throughout the paper we are dealing with finite sequences (or words), consequently we recall the following notation on sequences in $\F_p[T]$. Let $W=w_1,w_2,\ldots,w_n$ be such a finite sequence, then we set $\vert W\vert =n$ for the length of the word $W$. If we have two words $W_1$ and $W_2$, then $W_1,W_2$ denotes the word obtained by concatenation. Moreover, if $y\in \F_p^*$, then we define $y\cdot W$ as the following sequence$$y\cdot W = y w_1, y^{-1}w_2,\ldots, y^{(-1)^{n-1}}w_n.$$As usual, we denote by $[W]\in \F_p(T)$ the finite continued fraction $w_1+1/(w_2+1/(\dots ))$. In this formula the $w_i$, called the partial quotients, are non constant polynomials. Still, we will also use the same notation if the $w_i$ are constant and the resulting quantity is in $\F_p$. However in this last case, by writting $[w_1,w_2,\dots,w_n]$ we assume that this quantity is well defined in $\F_p$, i.e. $w_n\neq 0,[w_{n-1},w_n]\neq 0,\dots,[w_2,\dots,w_n]\neq 0$.
\newline We use the notation $\langle W\rangle$ for the continuant built from $W$. We denote by $W'$ (resp. $W''$) the word obtained from $W$ by removing the first (resp. last) letter of $W$. Hence, we recall that we have $[W]=\langle W\rangle/\langle W'\rangle$. We let $W^*=w_n,w_{n-1},\ldots,w_1$, be the word $W$ written in reverse order. We also have $[W^*]=\langle W\rangle/\langle W''\rangle$. It is also known that $[y\cdot W]=y[W]$.
\par If $\alpha \in \F(p)$ is an infinite continued fraction, $\alpha=[a_1,a_2,\dots,a_n,\dots]$, we set $x_n=\langle a_1,a_2,\dots,a_n\rangle$ and $y_n=\langle a_2,\dots,a_n\rangle$. In this way, we have $x_n/y_n=[a_1,a_2,\dots,a_n]$, with $x_1=a_1$, $y_1=1$ and by convention $x_0=1$, $y_0=0$. Recall that, if $\alpha_{n+1}=[a_{n+1},a_{n+2},\dots]$ is the tail of the expansion, we have $\alpha=(x_{n}\alpha_{n+1}+x_{n-1})/(y_{n}\alpha_{n+1}+y_{n-1})$, for $n\geq 1$.
\par As above $p$ is an odd prime and $k$ an integer with $1\leq k<p/2$. Linked to the previous polynomial mentioned above: $P_k(T)=(T^2-1)^k$, we introduce a second polynomial $Q_k$ in $\F_p[T]$. We define$$\omega_k=(-1)^k2k\prod_{1\leq i\leq k}(1-1/2i)\in \F_p^* \quad \text{ and }\quad Q_k=\omega_k^{-1}(A_1P_k-T^p).$$This pair $(P_k,Q_k)$ of polynomials was introduced in [L2]. The second polynomial can also be defined by (see [L2, p. 341])$$Q_k(T)=\int_{0}^{T}(x^2-1)^{k-1}dx=\sum_{0\leq i\leq k-1}(-1)^{k-1-i}\binom{k-1}{i}(2i+1)^{-1}T^{2i+1}.$$ Note that we also have $Q_k(1)=-\omega_k^{-1}$. We recall the following stated in [L2, p. 332]. 
\par {\emph{Let $l\geq 1$ be an integer and $(a_1,a_2,\dots,a_l)\in (\F_p[T])^l$, with $\deg(a_i)>0$ for $1\leq i\leq l$. Let $r=p^t$ with $t>0$ and $(P,Q)\in (\F_p[T])^2$ with $\deg(Q)<\deg(P)<r$. Then there exist a unique infinite continued fraction $\alpha=[a_1,a_2,\dots,a_l,\alpha_{l+1}]\in \F(p)^{+}$ satisfying $(*)\quad \alpha^r=P\alpha_{l+1}+Q$. This element $\alpha$ is the unique root in  $\F(p)^{+}$ of the following algebraic equation:
$$(**)\quad y_lX^{r+1}-x_lX^r+(Py_{l-1}-Qy_l)X-Px_{l-1}+Qx_l=0.$$}}
\par As above, let $\alpha$ be defined by the $l$-tuple, $(a_1,a_2,\dots,a_l)$, and the equality $(*)$. Then we set the following definitions.
\par {\emph{$\bullet \quad \alpha$ is a $P_k$-expansion of first kind if $r=p$ and there exists $(\epsilon_1,\epsilon_2)\in (\F_p^*)^2$ such that $(P,Q)=(\epsilon_1P_k,\epsilon_2Q_k)$.}}
\par {\emph{$\bullet \quad \alpha$ is a $P_k$-expansion of second kind if $r=p^2$ and there exists $(\epsilon_1,\epsilon_2)\in (\F_p^*)^2$ such that $(P,Q)=(\epsilon_1P_k^{p-1},\epsilon_2Q_k^p)$.}}
\par The importance of the pair $(P_k,Q_k)$ is due to different properties. The first and main one is the following (see [L2, p. 341], and note that $\omega_k$ was defined differentlty there than it is here below).
\par {\emph{Let $W_1$ be the finite word such that $P_k/Q_k=[W_1]$. Then we have 
$$W_1=v_{1}T,\dots,v_{i}T,\dots,v_{2k}T,$$where the numbers $v_i\in \F_p^*$ are defined by $v_1=2k-1$ and recursively, for $1\leq i\leq 2k-1$, by
$$v_{i+1}v_{i}=(2k-2i-1)(2k-2i+1)(i(2k-i))^{-1}.$$
Furthermore we have $W_1=-\omega_k^2\cdot W_1^*$.}}
\par This last equality implies a basic lemma which is the key tool in the study of the $P_k$-expansions. This lemma allows to transform a particular rational function into a continued fraction. In this way, starting from $(*)$, the partial quotients in a $P_k$-expansion of first kind have been obtained explicitly. This lemma is the following (see the origin in [L2, p. 343]).
\newline {\bf{Lemma 1. }}{\emph{ Let $A\in \F_p[T]$, $\delta \in \F_p^*$ and $X\in \F(p)$. Then we have $$A+\delta Q_kP_k^{-1}+X=[A, \delta^{-1}\cdot W_1,X'],$$ where
$$X'=X^{-1}P_k^{-2}+\omega_k^2\delta^{-1}Q_kP_k^{-1}.$$}}
We use this lemma to establish the continued fraction for the rational $P_k^{p-1}/Q_k^p$. This continued fraction will be fundamental to study $P_k$-expansions of second kind. We have the following proposition.
\newline {\bf{Proposition 2. }}{\emph{ Let $W_2$ be the finite word such that $P_k^{p-1}/Q_k^p=[W_2]$. Then $W_2$ is obtained from $W_1$ in the following way:
$$W_2=v_{1}A_{1},w_1\cdot W_1,v_{2}A_{1},w_2\cdot W_1,\dots,v_{2k-1}A_{1},w_{2k-1}\cdot W_1,v_{2k}A_{1},$$
where the numbers $w_{i}\in \F_p^*$ are defined by 
$$w_i^{-1}=-\omega_k[v_i,v_{i-1},\dots,v_1]\quad \text{ for }\quad 1\leq i\leq 2k-1.$$}}
Proof: We have $P_kQ_k^{-1}=[v_{1}T,\dots,v_{2k}T]$. Since $P_k(1)=0$ and $Q_k(1)\neq 0$, we obtain $[v_{1},\dots,v_{2k}]=0$ and $[v_{i},\dots,v_{1}] \in \F_p^*$ for $1\leq i\leq 2k-1$. We set $\alpha=[W_1]$ and $\beta=[W_2]$. We have$$\beta=(P_k/Q_k)^pP_k^{-1}=\alpha^pP_k^{-1}=[v_{1}T^p,\dots,v_{2k}T^p]P_k^{-1}=[v_{1}T^pP_k^{-1},P_k\alpha_2^p].$$Since we have $T^p=A_1P_k-\omega_kQ_k$, the last equality becomes$$\beta=v_{1}A_1-\omega_kv_{1}Q_kP_k^{-1}+P_k^{-1}\alpha_2^{-p}.$$Applying Lemma 1, with $\delta =-\omega_kv_{1}$ and $X=P_k^{-1}\alpha_2^{-p}$, we obtain$$\beta=[v_{1}A_1,w_1\cdot W_1,X'],$$ where $$X'=P_k^{-1}\alpha_2^{p}+\omega_k^2w_1Q_kP_k^{-1}.$$Since $\vert \alpha_2\vert =\vert T\vert$, we have $\vert X'\vert >1$. Consequently, we get $$b_1,\dots,b_{2k+1}=v_{1}A_1,w_1\cdot W_1\quad \text{ and }\quad X'=\beta_{2k+2}.$$Then we have$$\beta_{2k+2}=[v_{2}T^pP_k^{-1},P_k\alpha_3^{p}]+\omega_k^2w_1Q_kP_k^{-1}.$$Using $T^p=A_1P_k-\omega_kQ_k$ and applying the same lemma, with $\delta =-\omega_kv_{2}+\omega_k^2w_1=w_2^{-1}$ and $X=P_k^{-1}\alpha_3^{-p}$, we get $$\beta_{2k+2}=[v_{2}A_1,w_2\cdot W_1,X'],$$where$$X'=P_k^{-1}\alpha_3^{p}+\omega_k^2w_2Q_kP_k^{-1}.$$As above, we get the desired partial quotients, from the rank $2k+2$ to the rank $4k+2$, and also $X'=\beta_{4k+3}$. The process carries on until we get$$\beta_{4k^2-2k-1}=[v_{2k-1}A_1,w_{2k-1}\cdot W_1,X'],$$where $$\beta_{4k^2}=X'=P_k^{-1}\alpha_{2k}^{p}+\omega_k^2w_{2k-1}Q_kP_k^{-1}=v_{2k}T^pP_k^{-1}+\omega_k^2w_{2k-1}Q_kP_k^{-1}.$$Again, using $T^p=A_1P_k-\omega_kQ_k$, this becomes$$\beta_{4k^2}=v_{2k}A_1+\omega_kQ_kP_k^{-1}(\omega_kw_{2k-1}-v_{2k}).$$But we have$$\omega_kw_{2k-1}-v_{2k}=-[v_{2k-1},\dots,v_{1}]^{-1}-v_{2k}=-[v_{2k},v_{2k-1},\dots,v_{1}]=0.$$Hence $\beta_{4k^2}=v_{2k}A_1$ and the proof of the proposition is complete.\par We make a last remark on the continued fraction for $P_k^{p-1}/Q_k^p$. As for $W_1$, it can be seen that we also have $W_2=-\omega_k^2\cdot W_2^*$. It follows from this equality that the very same lemma as Lemma 1, holds in the second case when the pair $(P_k,Q_k)$ is replaced by the pair $(P_k^{p-1},Q_k^p)$ and $W_1$ is replaced by $W_2$. In this way, $P_k$-expansions of second kind could possibly be studied following the same path as in the first case ([L2] and [L4]). 
\vskip 0.5 cm
\noindent{\bf{3. Perfect $P_k$-expansions }}
\par If $p$ and $k$ are fixed, we recall that a $P_k$-expansion is defined by the $(l+2)$-tuple $(a_1,a_2,\dots,a_l,\epsilon_1,\epsilon_2)\in (\F_p[T])^l\times (\F_p^*)^2$. Once such a $(l+2)$-tuple is fixed, from the algebraic equation $(**)$, a computer can give the first partial quotients of the expansion, then it appears that some expansions are more ''regular'' than others. To be more precise, we use the following terminology. If $P\in \F_p[T]$, we say that $P$ is of type A (resp. of type B) if there exist $\lambda \in \F_p^*$ and $n\in \N$ such that $P=\lambda A_n$ (resp. $P=\lambda B_n$), where the polynomials $A_n$ (or $B_n$) belong to the sequences defined in the introduction. Then we say that a $P_k$-expansion of first kind (resp. of second kind) is perfect if every partial quotient is of type A (resp. of type B). It appears that a sufficient condition on the $(l+2)$-tuple $(a_1,a_2,\dots,a_l,\epsilon_1,\epsilon_2)$ can be given, in order to have a perfect $P_k$-expansion and then a description of the corresponding sequence of partial quotients. This is what is discussed in this section, distinguishing each of both cases.
\par We will use the following notation. For $n\geq 0$, if we have the word $w,w,\dots,w$ of length $n$, then we denote it shorly by $w^{[n]}$ with $w^{[0]}=\emptyset$. In the same way $W^{[n]}$ denotes the world $W,W,\dots,W$ where $W$ is repeated $n$ times and $W^{[0]}=\emptyset$. If we have a finite sequence $W=w_1,w_2,\ldots,w_n$ of polynomials of type A (or of type B), we can associate it to a finite sequence $I=i_1,i_2,\ldots,i_n$ of positive integers, such that $w_m=\lambda_m A_{i_m}$ (or  $w_m=\lambda_m B_{i_m}$), with $\lambda_m\in \F_p^*$ for $1\leq m\leq n$. Let $I_1$ (resp. $I_2$) denote the sequence of integers attached to the word $W_1$ (resp. $W_2$) introduced in the previous section. Then we have
$$I_1=0^{[2k]} \quad \text{ and }\quad I_2=1,(0^{[2k]},1)^{[2k-1]}.$$
Our aim is to describe the infinite sequence of integers associated to the infinite sequence of partial quotients for a perfect $P_k$-expansion. In each case, the sequences $I_1$ or $I_2$ are the stones from which this sequence is built. 
\par{\bf{A) Perfect $P_k$-expansions of first kind}} 
\par Let us consider a $P_k$-expansion of first kind. We have the following statement.
\newline {\it{ Assuming that the $(l+2)$-tuple $(a_1,\dots,a_l,\epsilon_1,\epsilon_2)$ satisfies an hypothesis $\mathcal{H}(1)$, then there exist a sequence $(\lambda_n)_{n\geq 1}$ in $\F_p^*$ and a sequence $(i(n))_{n\geq 1}$ in $\N$, such that
$$a_n=\lambda_n A_{i(n)}\quad \text{ for }\quad n\geq 1.\eqno{(1)}$$}} 
Note that the hypothesis $\mathcal{H}(1)$ will only be a sufficient condition to have $(1)$. This first case has been studied in previous works, and also in a more general setting (see [L2] and [L4]). A first task is to choose the first $l$ partial quotients $(a_1,\dots,a_l)$. Note that they must be of type A. In the case of the solution of $(Eq)$, for $p\equiv 1 \mod 3$, the first partial quotients are linear, hence we shall here consider them proportional to $A_0=T$. However, note that a more general situation could have been considered, as this was remarked in [L5]. Moreover, the extremal case $k=(p-1)/2$ conduces to perfect expansions having all partial quotients proportional to $T$. In a joint work with J.-Y. Yao, the second author, following a similar method, could obtain $P_k$-expansions of first kind having all partial quotients of degree 1, starting from $l$ partial quotients of degree 1, not necessarily proportional to $T$. 
\par Let us now describe our hypothesis $\mathcal{H}(1)$ concerning the $(l+2)$-tuple $(a_1,\dots,a_l,\epsilon_1,\epsilon_2)$. 
\newline {\it{We assume that there exists a $(l+2)$-tuple $(\lambda_1,\dots,\lambda_l,\epsilon_1,\epsilon_2)$ in $(\F_p^*)^{l+2}$ such that we have firstly $(a_1,\dots,a_l,)=(\lambda_1T,\dots,\lambda_lT)$ and secondly $$[\lambda_l,\lambda_{l-1},\dots,\lambda_1+\omega_k^{-1}\epsilon_2]=2k\epsilon_1\epsilon_2^{-1}.$$}} 
Note, as indicated above, that the existence of the square bracket in this last formula implies a special choice of the $(l+2)$-tuple. Indeed, one can check that there are exactly $(p-1)(p-2)^l$ such $(l+2)$-tuples in $(\F_p^*)^{l+2}$.
\newline In [L4, p. 256], in a larger context, it has been proved that the hypothesis $\mathcal{H}(1)$ implies that the $P_k$-expansion, defined by the corresponding $(l+2)$-tuple, is perfect, i.e. $(1)$ is satisfied. Both sequences $(\lambda_n)_{n\geq 1}$ and $(i(n))_{n\geq 1}$ have been described there. We do not give here indications on the sequence $(\lambda_n)_{n\geq 1}$ which is obtained by sophisticated recursive formulas from the initial $l$-tuple $(\lambda_1,\lambda_2,\dots,\lambda_l)$ and the pair $(\epsilon_1,\epsilon_2)$. Concerning the other sequence $(i(n))_{n\geq 1}$, we will describe it here simply with the above notation. Let  $(V_n)_{n\geq 0}$ be the sequence of finite words of integers defined recursively by
$$V_0=0\quad \text{ and }\quad V_n=n,V_0^{[2k]},V_1^{[2k]},\dots,V_{n-1}^{[2k]},\quad \text{ for }n\geq 1.\eqno{(2)}$$
Then the sequence $(i(n))_{n\geq 1}$ in $\N$ is given by the infinite word:
$$V_0^{[l]},V_1^{[l]},V_2^{[l]},\dots,V_n^{[l]},\dots \eqno{(3)}$$
\vskip 0.5 cm
\par{\bf{B) Perfect $P_k$-expansions of second kind}}
\par Let us consider a $P_k$-expansion of second kind. We have the following conjecture.
\newline {\it{ Assuming that the $(l+2)$-tuple $(a_1,\dots,a_l,\epsilon_1,\epsilon_2)$ satisfies an hypothesis $\mathcal{H}(2)$, then there exist a sequence $(\lambda_n)_{n\geq 1}$ in $\F_p^*$ and a sequence $(i(n))_{n\geq 1}$ in $\N$, such that$$a_n=\lambda_nB_{i(n)}\quad \text{ for }\quad n\geq 1.\eqno{(4)}$$}}  
Before describing the hypothesis $\mathcal{H}(2)$, we make a comment on the origin of this conjecture. We started from the solution of $(Eq)$ for $p=5, 11$ or $17$ and we obtained with a computer several thousands of partial quotients. This was not enough to guess the general pattern of this sequence of partial quotients. Inspired by the first case, we expected these particular expansions to belong to a much larger family. For small values of $p$, we knew that these expansions were generated in the way described above, i.e. were $P_k$-expansions of second kind (see [L5, p. 33]). By observing the first $l$ partial quotients and the pair $(\epsilon_1,\epsilon_2)$ in these three cases, we could guess the right form of the $(l+2)$-tuple $(a_1,\dots,a_l,\epsilon_1,\epsilon_2)$ in order to have a perfect expansion. It was only then, by considering a general $(l+2)$-tuple and letting the parameters $p,k$ and $l$ vary, that firstly we could state the above conjecture and secondly we could describe the pattern for these continued fractions.
\par  Let us now describe our hypothesis $\mathcal{H}(2)$ concerning the $(l+2)$-tuple $(a_1,\dots,a_l,\epsilon_1,\epsilon_2)$. 
\newline {\it{ Let $m\geq 1$ and $m$ integers $n_1,n_2\dots,n_m$ be such that
$$1<n_1<n_2<\dots<n_m \quad \text{ with }\quad n_{i+1}-n_i\geq 3 \quad \text{ for }\quad 1\leq i<m.$$
We set $l=n_m$ and we consider a $l$-tuple $(\lambda_1,\lambda_2,\dots,\lambda_l)\in (\F_p^*)^l$ such that 
$$[\lambda_1,\dots,\lambda_{n_1-1}]\neq 0\quad \text{ and }\quad [\lambda_{n_i+1},\dots,\lambda_{n_{i+1}-1}]=0 \quad \text{ for }\quad 1\leq i<m.$$
Then we assume that, for $1\leq n\leq l$ and $i=1\dots m$, we have $$n\neq n_i \Rightarrow a_n=\lambda_nB_0 \quad \text{ and }\quad n=n_i \Rightarrow a_n=\lambda_nB_1$$ and also $$\epsilon_2=-\omega_k[\lambda_1,\dots,\lambda_{n_1-1}]. $$}}
Note that the $m$-tuple $(\lambda_{n_1},\lambda_{n_2}\dots,\lambda_{n_m})\in (\F_p^*)^{m}$ as well as $\epsilon_1\in \F_p^*$ are chosen arbitrarily. Recall that, according to the remark made above, the existence of the square brackets in $\F_p$ implies restrictions on the choice of the remaining $\lambda_i$.
\par Now, we can describe conjecturally the sequence $(i(n))_{n\geq 1}$ under the hypothesis $\mathcal{H}(2)$. We define $$l_1=n_1-1 \quad \text{ and }\quad l_{i+1}=n_{i+1}-n_i-1 \quad \text{ for }\quad 1\leq i<m.$$ Let $V_0$ be the sequence of integers attached to the first $l$ partial quotients. According to the description made above, we can write
$$V_0=0^{[l_1]},1,0^{[l_2]},1,\dots,1,0^{[l_m]},1.\eqno{(5)}$$
We introduce the following sequence $(J_n)_{n\geq 1}$ of words, defined recursively by $J_1=0^{[2k]},1$ and 
$$J_{n+1}=(2n,2n-1,J_n^{[2k-1]})^{[2k-1]},2n,2n+1\quad \text{ for }n\geq 1.\eqno{(6)}$$
Next, for $1\leq i\leq m$ and $n\geq 1$, we define 
$$V_{n,i}=(2n,2n-1,J_n^{[2k-1]})^{[l_i-1]},2n,2n+1. \eqno{(7)}$$
Then the sequence $(i(n))_{n\geq 1}$ in $\N$ is given by the infinite word:
$$V_0,V_1,V_2,\dots,V_n,\dots \quad \text{ where }\quad V_n=V_{n,1},V_{n,2},\dots,V_{n,m}\quad \text{ for }n\geq 1.\eqno{(8)}$$
This conjecture on $P_k$-expansions of second kind gives each partial quotient up to a multiplicative constant in $\F_p^*$. Note that this gives the degree of each partial quotient. We will see in the last section that the conjecture is conforted by considerations on the growth of this sequence of degrees. Finally, let us repeat that we have not tried to describe, even conjecturally, the sequence of constants in (4), which depends on the $(l+1)$-tuple $(\lambda_1,\dots,\lambda_l,\epsilon_1)$. 
\vskip 0.5 cm
\noindent{\bf{4. Link between $P_k$-expansions and the solution of $(Eq)$ }}
\par Let $p>3$ be a prime number. We set $k=(p-1)/3$ if $p\equiv 1 \mod 3$, and otherwise $k=(p+1)/3$. We set $(r,P,Q)=(p,P_k,Q_k)$ if $p\equiv 1 \mod 3$, and otherwise $(r,P,Q)=(p^2,P_k^{p-1},Q_k^p)$. We set $l=(p-1)/2$ if $p\equiv 1 \mod 3$, and otherwise $l=(p+1)^2/3$. We set ${\bf{\Lambda}}(p)=(a_1,\dots,a_l,\epsilon_1,\epsilon_2)$, where $a_i\in \F_p[T]$, with $\deg(a_i)>0$ for $1\leq i\leq l$, and $(\epsilon_1,\epsilon_2)\in (\F_p^*)^2$. Then we define in $\F_p[T][X]$ the following polynomial:
$$H({\bf{\Lambda}}(p);X)=y_lX^{r+1}-x_lX^r+(\epsilon_1Py_{l-1}-\epsilon_2Qy_l)X-\epsilon_1Px_{l-1}+\epsilon_2Qx_l.$$
Recall that $x_l,y_l,x_{l-1}$ and $y_{l-1}$ are the continuants built from the $l$-tuple $(a_1,\dots,a_l)$, as stated at the beginning of Section 2. Let $A\in \F_p[T][X]$ be the polynomial defined by $A(X)=(9/32)X^4-TX^3+X^2-8/27$.
\par Our goal is to show that there exists a particular $(l+2)$-tuple, denoted by ${\bf{\Lambda}_0}(p)$, such that $A(X)$ divides $H({\bf{\Lambda}_0}(p);X)$ in the ring of polynomials with coefficients in $\F_p[T]$. The solution $\alpha(p)$ of $(Eq)$ satisfies $A(\alpha(p))=0$, consequently we will have $H({\bf{\Lambda_0}}(p);\alpha(p))=0$, so that $\alpha(p)$ will be a $P_k$-expansion of first kind if $p\equiv 1 \mod 3$ and of second kind if $p\equiv 2 \mod 3$, according to what is stated in Section 2. Moreover, this particular ${\bf{\Lambda}_0}(p)$ will be chosen such that $\mathcal{H}(1)$ is satisfied if $p\equiv 1 \mod 3$, and otherwise $\mathcal{H}(2)$ is satisfied. In the first case, $\mathcal{H}(1)$ being satisfied, this will prove that $\alpha(p)$ is a perfect $P_k$-expansion of first kind, implying $(I)$. In the second case, $\mathcal{H}(2)$ being satisfied, this will imply the conjecture stating that $\alpha(p)$ is a perfect $P_k$-expansion of second kind, and consequently satisfies $(II)$. However this argument is based on the division of $H({\bf{\Lambda}_0}(p);.)$ by $A$ and this division is established by straightforward computations with the help of a computer. It follows that in both cases (proof or conjecture) the result holds with a bound on the prime number $p$, eventhough there is no reason to doubt that the same is true for all primes $p$. 
\par  Now we shall describe the $(l+2)$-tuple ${\bf{\Lambda}_0}(p)$, in order to obtain the divisibility property mentioned above. We set ${\bf{\Lambda}_0}(p)=(B(p),\epsilon_1,\epsilon_2)$. We shall first describe the word $B(p)=a_1,\dots,a_l$. By observations on the continued fraction expansion of the solution of $(Eq)$, for small values of $p$, we could guess the form of $B(p)$. Let us consider the words $W_1$ and $W_2$ defined in Section 2. For $p\equiv 1 \mod 3$ and $k=(p-1)/3$, we let $\widehat{W_1}$ be the word formed by the last letters of $W_1$, as follows:
$$\widehat{W_1}=v_{k/2+1}T,v_{k/2+2}T,\dots,v_{2k}T \quad \text{ with }\quad \vert \widehat{W_1}\vert =3k/2=(p-1)/2.$$
In the same way, for $p\equiv 2 \mod 3$ and $k=(p+1)/3$, according to Proposition 2, we have $W_2=b_1,\dots,b_{4k^2}$. Then we let $\widehat{W_2}$ be defined by
$$\widehat{W_2}=b_{k^2+1},b_{k^2+2},\dots,b_{4k^2} \quad \text{ with }\quad \vert \widehat{W_2}\vert =3k^2=(p+1)^2/3.$$
Note that $\vert B(p) \vert=l$, and consequently, in each of the corresponding case, we have $\vert B(p) \vert=\vert \hat{W_1}\vert$ or $\vert B(p) \vert=\vert \hat{W_2}\vert$. We noted that, for small values of $p$, in each of both cases, the word $B(p)$ is equal to $\widehat{W_1}$ or $\widehat{W_2}$ up to a small correction. Consequently, for all $p$, with $\epsilon$ in $\F_p^*$, we define $B(p)$ by 
$$B(p)=\epsilon\cdot \widehat{W_1} \quad \text{ if }\quad p\equiv 1 \mod 3\eqno{(9)}$$
and 
$$B(p)=\epsilon\cdot \widehat{W_2} \quad \text{ if }\quad p\equiv 2 \mod 3.\eqno{(10)}$$
Recall that, for all $p$, the solution $\alpha(p)$ of $(Eq)$ is such that $a_1=(32/9)T$. By identification from $(9)$ and $(10)$, in each case $p \equiv 1$ or $2\mod 3$, we see that $\epsilon$ must satisfy
$$\epsilon v_{k/2+1}T=(32/9)T\quad \text{ or }\quad \epsilon b_{k^2+1}=(32/9)T.$$
From Proposition 2, one can check that we have $b_{k^2+1}=w_{k/2}^{(-1)^{k/2}}v_{3k/2+1}T$. Consequently, $\epsilon \in \F_p^*$ will be defined by
$$\epsilon=32/(9v_{k/2+1}) \quad \text{ if }\quad p\equiv 1 \mod 3\eqno{(11)}$$
and 
$$\epsilon=32/(9w_{k/2}^{(-1)^{k/2}}v_{3k/2+1}) \quad \text{ if }\quad p\equiv 2 \mod 3.\eqno{(12)}$$
To fully determine ${\bf{\Lambda}_0}(p)$, we must now describe the pair $(\epsilon_1,\epsilon_2)$ in $(\F_p^*)^2$.  This will be achieved in two steps. First, we introduce an argument based on experimental observation.
\par We set $U=\epsilon_1Py_{l-1}-\epsilon_2Qy_l$ and $V=\epsilon_2Qx_l-\epsilon_1Px_{l-1}$. Consequently the polynomial $H$ considered above is $H=y_lX^{r+1}-x_lX^r+UX+V$. We have observed, for the first values of $p$ (in both cases), that the polynomials $U$ and $V$ have a low degree in the indeterminate $T$. This will imply a useful relation between $\epsilon_1$ and $\epsilon_2$. Let $a\in \F_p[T]$ be defined by $a=T$ if $p\equiv 1 \mod 3$ and $a=A_1$ if $p\equiv 2 \mod 3$. From the continued fraction expansion of $P/Q$, for $p \equiv 1$ or $2\mod 3$, we get
$$P=v_1aQ+R \quad \text{ where }\quad R\in \F_p[T] \quad \text{ and }\quad \vert R
\vert <\vert Q\vert.\eqno{(13)}$$
From $(9)$ and $(10)$, and from the continued fraction expansion of $P/Q$, again for $p \equiv 1$ or $2\mod 3$, we have
$$a_l=v_{2k}a\epsilon^{(-1)^{l+1}}.\eqno{(14)}$$
We have $x_l=a_lx_{l-1}+x_{l-2}$ and $y_l=a_ly_{l-1}+y_{l-2}$. Combining these equalities with $(13)$ and $(14)$, we obtain
$$U=(\epsilon_1v_1-\epsilon_2v_{2k}\epsilon^{(-1)^{l+1}})aQy_{l-1}+\epsilon_1Ry_{l-1}-\epsilon_2Qy_{l-2}\eqno{(15)}$$
and
$$V=(\epsilon_2v_{2k}\epsilon^{(-1)^{l+1}}-\epsilon_1v_1)aQx_{l-1}+\epsilon_2Qx_{l-2}-\epsilon_1Rx_{l-1}.\eqno{(16)}$$
Since $\vert R\vert <\vert P\vert$ and $\vert y_{l-2}\vert <\vert y_l\vert$, we get easily $\vert \epsilon_1Ry_{l-1}-\epsilon_2Qy_{l-2}\vert <\vert U\vert$. In the same way we have $\vert \epsilon_2Qx_{l-2}-\epsilon_1Rx_{l-1}\vert <\vert V\vert$. In order to have $\vert U\vert$ and $\vert V\vert$ small, according to $(15)$ and $(16)$, we must assume the following 
$$\epsilon_2v_{2k}\epsilon^{(-1)^{l+1}}=\epsilon_1v_1.\eqno{(17)}$$
Finally, since $v_1=-\omega_k^2v_{2k}$, , by $(17)$ for $p \equiv 1$ or $2\mod 3$, we assume that
$$\epsilon_1=-\epsilon_2\omega_k^{-2}\epsilon^{(-1)^{l+1}}.\eqno{(18)}$$
\par In a last step, according to the case considered, we shall make sure thare $\mathcal{H}(1)$ or $\mathcal{H}(2)$ is satisfied. 
\newline Let us consider the first case: $p \equiv 1 \mod 3$. According to the form of $B(p)$, we have $a_i=\lambda_i T=\epsilon^{(-1)^{i+1}}v_{k/2+i}T$ for $1\leq i\leq l$. To have $\mathcal{H}(1)$ satisfied, we need have
$$[\lambda_l,\lambda_{l-1},\dots,\lambda_1+\omega_k^{-1}\epsilon_2]=2k\epsilon_1\epsilon_2^{-1}. $$ 
This formula can be inversed and it is equivalent to
$$[\lambda_1,\lambda_2,\dots,\lambda_l-2k\epsilon_1\epsilon_2^{-1}]=-\omega_k^{-1}\epsilon_2.   \eqno{(19)}$$
But, using $(9)$ and recalling that $[\epsilon \cdot W]=\epsilon [W]$, we obtain
  $$[\lambda_1,\lambda_2,\dots,\lambda_l-2k\epsilon_1\epsilon_2^{-1}]=\epsilon[v_{k/2+1},\dots,v_{2k}-2k\epsilon_1\epsilon_2^{-1}\epsilon^{(-1)^{l}}].\eqno{(20)}$$
By $(18)$ and $v_1=2k-1=-\omega_k^2v_{2k}$, with $k=-1/3$ in $\F_p^*$, we can write
$$v_{2k}-2k\epsilon_1\epsilon_2^{-1}\epsilon^{(-1)^{l}}=v_{2k}+2k\omega_k^{-2}=-v_{2k}/(2k-1)=3v_{2k}/5. \eqno{(21)}$$
Combining $(19)$, $(20)$ and $(21)$, we see that $\mathcal{H}(1)$ is satisfied if 
$$\epsilon_2=-\epsilon\omega_k[v_{k/2+1},v_{k/2+2},\dots,v_{2k-1},3v_{2k}/5].\eqno{(22)}$$
Let us now consider the second case: $p \equiv 2 \mod 3$. We first check that $B(p)$ has the form required in $\mathcal{H}(2)$. According to Proposition 2 and the definition of $\widehat{W_2}$, we have 
$$\widehat{W_2}=W_3,v_{k/2+1}B_{1},w_{k/2+1}\cdot W_1,v_{k/2+2}B_{1},\dots,v_{2k-1}B_{1},w_{2k-1}\cdot W_1,v_{2k}B_{1},$$ 
where
$$W_3=w_{k/2}^{(-1)^{k/2}}v_{3k/2+1}T,w_{k/2}^{(-1)^{k/2+1}}v_{3k/2+2}T,\dots,w_{k/2}^{-1}v_{2k}T.$$
Since $B(p)=\epsilon \cdot \widehat{W_2}$, we see that $a_i=\lambda_iB_0$ or $a_i=\lambda_iB_1$, and $B(p)$ has the right form. More precisely, with the notation introduced in section 3, we have $m=3k/2=(p+1)/2$ and, for $1\leq i\leq m$, $n_i=i(2k+1)-3k/2$. Also $\vert W_3\vert=k/2=l_1$ and, for $1< i\leq m$, $l_i=2k$. Moreover the following equalities hold
$$[\lambda_1T,\dots,\lambda_{n_1-1}T]=[\epsilon \cdot W_3]\eqno{(23)}$$
and, for $1\leq i<m$,
$$[\lambda_{n_i+1}T,\dots,\lambda_{n_{i+1}-1}T]=[\epsilon^{\pm 1}\cdot (w_{k/2+i}\cdot W_1)].\eqno{(24)}$$
Since $[W_1]_{T=1}=0$, for $1\leq i<m$, we observe that $(24)$ implies the required condition
$$[\lambda_{n_i+1},\dots,\lambda_{n_{i+1}-1}]=\epsilon^{\pm 1}w_{k/2+i}[W_1]_{T=1}=0.$$
Finally, by $(23)$, we have
$$[\lambda_1,\dots,\lambda_{n_1-1}]=\epsilon [W_3]_{T=1}=\epsilon w_{k/2}^{(-1)^{k/2}}[v_{3k/2+1},\dots,v_{2k}].\eqno{(25)}$$
Therefore $\mathcal{H}(2)$ is fully satisfied if we fix
$$\epsilon_2=-\epsilon\omega_k w_{k/2}^{(-1)^{k/2}}[v_{3k/2+1},v_{3k/2+2},\dots,v_{2k}].\eqno{(26)}$$
\par We have completed the description of ${\bf{\Lambda}_0}(p)$ in both cases. If $p \equiv 1 \mod 3$, we apply $(9)$, $(11)$, $(18)$ and $(22)$. If $p \equiv 2 \mod 3$, we apply $(10)$, $(12)$, $(18)$ and $(26)$. From ${\bf{\Lambda}_0}(p)$, we can build the coefficients in $\F_p[T]$ of $H$ and then, for a given prime $p$, it is a routine work to check by computer that the polynomial $A$ divides the polynomial $H$ in $\F_p[T][X]$. We have done so for the first primes in both cases.
\par We make a last comment on the polynomial $H$. In the first case, in [L5, p. 31], a different form of $H$ has been given. The polynomials $U$ and $V$ were expicited using continuants involving $v_1T,v_2T,\dots,v_{k/2}T$. As a consequence, in this case, we have $\vert U\vert=\vert T\vert^{k/2}$ and $\vert V\vert=\vert T\vert^{k/2-1}$. In the second case, the same approach is possible and this leads to $\vert U\vert=\vert T\vert^{3k^2/2-k}$ and $\vert V\vert=\vert T\vert^{3k^2/2-k-1}$. 
\vskip 0.5 cm
\noindent{\bf{5. Irrationality measure}}
\par Mahler [M], by an adaptation of an old and famous theorem on algebraic real numbers, due to Liouville, could prove the following theorem.
\newline {\it{Let $\alpha \in \F(p)$ be algebraic of degree $n>1$ over $\F_p(T)$. Then there exists a positive real constant $C$ such that
$$ \vert \alpha -P/Q \vert \geq C \vert Q \vert^{-n} \quad \text{ for all }\quad P/Q \in \F_p(T).$$}}
The irrationality measure of an irrational power series $\alpha$ is defined by
$$\nu(\alpha)=-\limsup_{\vert Q \vert \to \infty}\log(\vert \alpha -P/Q \vert)/\log(\vert Q \vert).$$
By Roth's theorem, the irrationality measure of an irrational algebraic real number is 2. For power series over a finite field there is no analogue of Roth's theorem. However, according to Mahler's theorem, we have $\nu(\alpha)\in [2,n]$ if $\alpha$ is algebraic of degree $n>1$ over $\F_p(T)$. 
\par The irrationality measure is directly related to the growth of the sequence of the degrees of the partial quotients $a_n$ in the continued fraction expansion of $\alpha$. Indeed we have (see [L1, p. 214])
 $$\nu(\alpha)=2+\nu_0(\alpha)=2+\limsup_n(\deg(a_{n+1})/\sum_{1\leq i\leq n}\deg(a_i)).\eqno{(27)}$$
This formula allows to compute $\nu(\alpha)$ if the continued fraction is explicitly known. Note that if $\alpha$ is a $P_k$-expansion then we have $\nu_0(\alpha)\in [0,r-1]$ and  moreover if $\alpha(p)$ is the solution of $(Eq)$ we have $\nu_0(\alpha(p))\in [0,2]$. We shall compute the irrationality measure for perfect $P_k$-expansions and for the solution of $(Eq)$.
\par We need the following notation. Let $W$ be a finite sequence $W=w_1,w_2,\ldots,w_n$ of polynomials of type A (or of type B), associated as above to the finite sequence $I=i_1,i_2,\ldots,i_n$ of positive integers, such that $w_m=\lambda_m A_{i_m}$ (or  $w_m=\lambda_m B_{i_m}$), where $\lambda_m\in \F_p^*$ for $1\leq m\leq n$. We let $D(W)$ be the sum $\sum_{1\leq i\leq n}\deg(w_i)$ and, by extension, we define $D(I)=D(W)$. 
\par $\bullet$ Perfect $P_k$-expansion of first kind.
\newline Here, for $n\geq 1$,  we have $a_n=\lambda_nA_{i(n)}$ and consequently $\deg(a_n)=u_{i(n)}$. From the recurrent definition of $A_n$, we have $u_n=(p^n(p-2k-1)+2k)/(p-1)$ for $n\geq 0$. To compute the limit in $(27)$, we need to observe the first occurence of $u_{n}$ in the sequence of the degrees of the partial quotients, then compute the sum of all the degrees appearing before this term in this sequence. According to the description of $(i(n))_{n\geq 1}$ given by $(2)$ and $(3)$ in Section 3, we see that $n$ appears for the first time at the begining of $V_n$. Hence, by $(27)$, we have 
$$\nu_0(\alpha)=\lim_n (u_n/D(V_0^{[l]},V_1^{[l]},V_2^{[l]},\dots,V_{n-1}^{[l]})).\eqno{(28)}$$  
Let us compute $D(V_n)$. We have $D(V_0)=u_0=1$ and also by $(2)$
$$D(V_n)=D(n,V_0^{[2k]},V_1^{[2k]},\dots,V_{n-1}^{[2k]})=u_n+2k\sum_{1\leq i\leq n-1}D(V_i).\eqno{(29)}$$
By induction, from $(29)$,x it is elementary to verify that $D(V_n)=p^n$ holds, for $n\geq 0$. Consequently $(28)$ implies
$$\nu_0(\alpha)=\lim_n (u_n/l\sum_{1\leq i\leq n-1}p^i)=\lim_n (p^n(p-2k-1)+2k)/(l(p^n-1)).$$
Finally, we obtain
$$\nu_0(\alpha)=(p-2k-1)/l.\eqno{(30)}$$
Note that $\nu_0(\alpha)=0$ if $2k=p-1$. In this case, we have $a_n=\lambda_nT$, for $n\geq 1$. However the sequence $(\lambda_n)_{n\geq 1}$ is generally not ultimately periodic and therefore $\alpha$ is not quadratic. See [LY], for a general study on algebraic continued fractions with partial quotients of degree 1.
\par $\bullet$ Perfect $P_k$-expansion of second kind.
\newline Here, for $n\geq 1$,  we have $a_n=\lambda_nB_{i(n)}$ and consequently $\deg(a_n)=v_{i(n)}$. We have $v_n=(p^n(p-2k+1)+2k(-1)^n)/(p+1)$ for $n\geq 0$, again from the recurrent definition of $B_n$. The description of the sequence $(i(n))_{n\geq 1}$ is given by the formulas $(5)-(8)$ in Section 3. We observe that $2n$ appears for the first time at the beginning of $V_n$, whereas $2n+1$ apears for the first time at the end of $V_{n,1}$. As above, to use $(27)$, we need to compute the sum of the degrees before these occurences. We set $V_{n,1}=W_{n,1},2n+1$ and we have $D(W_{n,1})=D(V_{n,1})-v_{2n+1}$. We set $r_n=D(V_0,V_1,\dots,V_{n-1})$. We have $D(V_0,V_1,\dots,V_{n-1},W_{n,1})=r_n+D(V_{n,1})-v_{2n+1}$. Then we define
$$s_n=v_{2n}/r_n\quad \text{and} \quad t_n=v_{2n+1}/(r_n+D(V_{n,1})-v_{2n+1}).\eqno{(31)}$$
Hence, by $(27)$, we have
$$\nu_0(\alpha)=\max(\lim_n s_n,\lim_n t_n). \eqno{(32)}$$
We need to compute $D(V_{n,i})$ and $D(V_{n})$. First, let us compute $D(J_n)$. We have $D(J_1)=D(0^{[2k]},1)=2kv_0+v_1=2k+p-2k=p$. For $n\geq 1$, by $(6)$, we get
$$D(J_{n+1})=(2k-1)(v_{2n}+v_{2n-1}+(2k-1)D(J_n))+v_{2n}+v_{2n+1}.\eqno{(33)}$$   
For $n\geq 1$, we observe that $v_n+v_{n-1}=p^{n-1}(p-2k+1)$. Hence $(33)$ can be written as
$$D(J_{n+1})=(2k-1)(p^{2n-1}(p-2k+1)+(2k-1)D(J_n))+p^{2n}(p-2k+1).$$ 
By induction, it is elementary to check that this formula implies
$$D(J_n)=p^{2n-1}\quad \text{for} \quad n\geq 1. \eqno{(34)}$$  
Now we use $(7)$ to compute $D(V_{n,i})$ for $n\geq 1$ and $1\leq i\leq m$. We have
$$D(V_{n,i})=(l_i-1)(v_{2n}+v_{2n-1}+(2k-1)D(J_n))+v_{2n}+v_{2n+1}.$$ 
By $(34)$ and the formula for $v_n+v_{n-1}$, this last formula gives
$$D(V_{n,i})=(l_i-1)p^{2n}+p^{2n}(p-2k+1)=p^{2n}(p-2k+l_i).\eqno{(35)}$$  
 Consequently, by $(8)$ and $(35)$, by sommation, with $\sum_{i=1}^{ m} l_i=l-m$, we get
$$D(V_{n})=p^{2n}\sum_{1\leq i\leq m}(p-2k+l_i)=p^{2n}(m(p-2k-1)+l).\eqno{(36)}$$
We set $\nu_1=m(p-2k-1)+l$. By $(5)$, $D(V_0)=\sum_{i=1}^{ m} l_i+m(p-2k)$ and consequently we have $D(V_0)=\nu_1$. Again by sommation, $(36)$ implies
$$r_n=\sum_{0\leq i\leq n-1}D(V_{i})=\nu_1(p^{2n}-1)/(p^2-1).\eqno{(37)}$$
Thus we have $s_n=v_{2n}/r_n=(p-1)(p^{2n}(p-2k+1)+2k)/(\nu_1(p^{2n}-1))$,
which implies
$$\lim_n s_n=(p-1)(p-2k+1)/\nu_1.\eqno{(38)}$$
We have $1/t_n=(1/s_n)(v_{2n}/v_{2n+1})+D(V_{n,1})/v_{2n+1}-1$. From $(35)$, we get
$$\lim_n (D(V_{n,1})/v_{2n+1})=(p-2k+l_1)(p+1)/(p(p-2k+1)).\eqno{(39)}$$
From $(38)$, since $\lim_n (v_{2n}/v_{2n+1})=1/p$, we also get 
$$\lim_n ((1/s_n)(v_{2n}/v_{2n+1}))=\nu_1/(p(p-1)(p-2k+1)).\eqno{(40)}$$
Combining $(39)$ and $(40)$, we obtain
$$\lim_n (1/t_n)=(\nu_1+(p-2k+l_1)(p^2-1))/(p(p-1)(p-2k+1))-1.$$
 Finally, this implies
$$\lim_n t_n=p(p-1)(p-2k+1)/\nu_2,\eqno{(41)}$$
where
$$\nu_2=\nu_1+(p-2k+l_1)(p^2-1)-p(p-1)(p-2k+1)=\nu_1+l_1(p^2-1)-2k(p-1).$$
Combining $(32)$, $(38)$ and $(41)$, we obtain
$$\nu_0(\alpha)=(p-1)(p-2k+1)\max (1/\nu_1,p/\nu_2).\eqno{(42)}$$
A simple calculation shows that $p-2k+1\leq \nu_1$ and $p(p-2k+1)\leq \nu_2$. Moreover we have $\nu_1=\nu_2/p=p-2k+1$ if and only if $m=1$ and $l=2$. Consequently, we obtain $\nu_0(\alpha)=p-1$ if $m=1$ and $l=2$ and $\nu_0(\alpha)<p-1$ otherwise. Note that $\nu_0(\alpha)$ is far from the admitted upper bound $p^2-1$.
\newline The extremal case, corresponding to  $m=1$ and $l=2$, is noteworthy. In this case, we have $l_1=1$, $V_0=0,1$ and $V_n=V_{n,1}=2n,2n+1$, this shows that the infinite word describing the sequence $(i(n))_{n\geq 1}$ is simply $\N$. We have $\epsilon_2=-\omega_k\lambda_1$, and the continued fraction is defined by the triple $(\lambda_1,\lambda_2,\epsilon_1)$ in $(\F_p^*)^3$. There exists a sequence $(\lambda_n)_{n\geq 1}$ in $\F_p^*$ such that $a_n=\lambda_nB_{n-1}$, for $n\geq 1$. Amazingly enough, this sequence is simply defined by $\lambda_{2n+1}=\epsilon_1^{-n}\lambda_1$ and $\lambda_{2n+2}=\epsilon_1^{n}\lambda_2$, for $n\geq 0$.
\par $\bullet$ Solution of $(Eq)$ if $p\equiv 1 \mod 3$.
\newline Here $\alpha(p)$ has been proved (with a limitation on $p$) to be a perfect $P_k$-expansion of first kind, where $k=(p-1)/3$ and $l=(p-1)/2$. Applying $(30)$, we get $\nu(\alpha(p))=2+\nu_0(\alpha(p))=8/3.$ 
\par $\bullet$ Solution of $(Eq)$ if $p\equiv 2 \mod 3$.
\newline According to Section 4, here we conjecture that $\alpha(p)$ is a perfect $P_k$-expansion of second kind, where $k=(p+1)/3$, $l=(p+1)^2/3$, $m=(p+1)/2$ and $l_1=k/2$. A direct computation shows that $(p-1)(p-2k+1)=2\nu_1/3$, whereas $p(p-1)(p-2k+1)=2\nu_2$. Consequently, applying $(42)$, we have  
$\nu(\alpha(p))=2+\nu_0(\alpha(p))=4$. According to Mahler's theorem, $\alpha(p)$ being algebraic of degree 4, $\nu(\alpha(p))$ has the maximal possible theoretical value.
\vskip 0.5 cm
\noindent{\centerline{\bf{Acknowledgment}}}
\par We thank warmly Bill Allombert for his advices in programming and his help to insert the figures in this note.

\vskip 1 cm
\begin{tabular}{ll}Khalil AYADI\\D\'epartement de Math\'ematiques\\Facult\'e des Sciences de Sfax\\Sfax 3018, Tunisie\\E-mail: ayedikhalil@yahoo.fr\\
\\Alain LASJAUNIAS\\Institut de Math\'ematiques de Bordeaux  CNRS-UMR 5251
\\Universit\'e de Bordeaux \\Talence 33405, France \\E-mail: Alain.Lasjaunias@math.u-bordeaux1.fr\\\end{tabular}

\end{document}